\documentclass[11pt]{article}%[10pt]{article}%\documentclass
\usepackage{amsmath,amssymb,latexsym,epsf,mathabx}
\usepackage[all]{xy}
\usepackage{mathrsfs}
\usepackage{soul}
\usepackage{color}

\begin{document}

%% ----------------------------------------------------------------------
\newcommand{\nc}{\newcommand}
\newtheorem{Th}{Theorem}[section]%[section]%
\newtheorem{Def}[Th]{Definition}
\newtheorem{Lem}[Th]{Lemma}
\newtheorem{Con}[Th]{Condition}
\newtheorem{Pro}[Th]{Proposition}
\newtheorem{Cor}[Th]{Corollary}
\newtheorem{Rem}[Th]{Remark}
\newtheorem{Exm}[Th]{Example}
\newtheorem{Sc}[Th]{}
\def\Pf#1{{\noindent\bf Proof}.\setcounter{equation}{0}}
\def\bskip#1{{ \vskip 20pt }\setcounter{equation}{0}}
\def\sskip#1{{ \vskip 5pt }\setcounter{equation}{0}}
\def\mskip#1{{ \vskip 10pt }\setcounter{equation}{0}}
\def\bg#1{\begin{#1}\setcounter{equation}{0}}
\def\ed#1{\end{#1}\setcounter{equation}{0}}

\soulregister\cite7 % 针对\cite命令
\soulregister\citep7 % 针对\citep命令
\soulregister\citet7 % 针对\citet命令
\soulregister\ref7

%%%%%%%%%%%%%%%%%%%%%%%%%%%%%%%%%%%%%%%%%%%%%%%%%%%%%%%%%%%%%%%%%%%%%%%%%%%%%%%%%
%**************************标题、摘要、分类号、关键字**************************

\title{\bf  Gorenstein syzygy objects in extriangulated categories
\thanks{Supported by the National Natural Science Foundation of China (Grants No.11801004) and the Startup Foundation for Introducing Talent of AHPU (Grant No.2017YQQ016 and 2020YQQ067). }
}
%\footnotetext{
%E-mail:~zhangpy@ahpu.edu.cn}
\smallskip
\author{ Peiyu Zhang\thanks{Correspondent author}, Ming Chen and Dajun Liu\\ %$^\text{a}$\\
\footnotesize ~E-mail:~zhangpy@ahpu.edu.cn, yysxcm@126.com, liudajun@ahpu.edu.cn; \\
\footnotesize School of Mathematics and Physics,  Anhui Polytechnic University, Wuhu, China.
\\ Jiaqun Wei\\ %$^\text{b,c}$\\
\footnotesize ~E-mail:~weijiaqun@njnu.edu.cn;\\
\footnotesize  School of Mathematics Science, Nanjing Normal University, Nanjing, China.
}

\date{}
\maketitle
\baselineskip 15pt%16pt%14pt%15.5pt%\baselineskip  25.5pt %
%%%%%%%\hskip 18pt
%
% Abstract ------------------------------------------------------
%
\begin{abstract}
\vskip 10pt%
We give the definition of Gorenstein syzygy objects in extriangulated categories and obtain a characterization.
For a recollement $R(\mathcal{A'}, ~\mathcal{A}, ~\mathcal{A''})$ of extriangulated categories, we mainly show that
Gorenstein syzygy objects in $\mathcal{A'}$ and $\mathcal{A''}$ induce Gorenstein syzygy objects in $\mathcal{A}$ under certain conditions,
and prove that Gorenstein syzygy objects in $\mathcal{A}$ induce Gorenstein syzygy objects in $\mathcal{A'}$ and $\mathcal{A''}$ under certain conditions,

\mskip\


%\sskip\


\end{abstract}
%\smallskip
%
\vskip 30pt

\section{Introduction}%}
%{\noindent \Large\bf Introduction}
Exact categories and triangulated categories are two fundamental structures in algebra.
Extriangulated categories were introduced by Nakaoka and Palu as a simultaneous generalization of exact categories and triangulated categories in \cite{NP}.
And the authors given an extriangulated categories which are neither exact categories nor triangulated categories, see \cite{NP,ZZ}.
Many scholars have done a lot of research on extriangulated categories, for instance \cite{GMT,HZ,HZZ,HZZ1,LN,WWZ,ZZT} and so on.

Recollements of triangulated categories were introduced by Be\u{1}linson, Bernstein and Deligne \cite{BBD} in connection with derived categories
of sheaves on topological spaces, the research of recollements refer to \cite{CJM,FP,GJ,MH,MH1,PC,LYN,PC1} and so on. The authors give a simultaneous
generalization of recollements of triangulated categories and Abealian categories, which call recollements of extriangulated categories \cite{WWZ},
prove the result as \cite{CJM}. Since the extriangulated category is a generalization of both triangulated and exact categories, so one
may ask if we can do a similar job as triangulated and exact categories. In this paper, we mainly consider Gorenstein syzygy modules \cite{HH} in extriangulated category,
which we called Gorenstein syzygy objects and the gluing Gorenstein syzygy objects in a recollement under certain conditions.

We now outline the results of the paper. In section 2, we first recalled some definitions and properties of extriangulated categories. In section 3,
we mainly give the definition of Gorenstein syzygy objects in extriangulated categories and a characterization.
And we can obtain a Equivalent characterization about the Gorenstein projective dimension. In section 4, we mainly show that
Gorenstein syzygy objects in $\mathcal{A'}$ and $\mathcal{A''}$ induce Gorenstein syzygy objects in $\mathcal{A}$ under certain conditions,
and prove that Gorenstein syzygy objects in $\mathcal{A}$ induce Gorenstein syzygy objects in $\mathcal{A'}$ and $\mathcal{A''}$ under certain conditions,

\section{Preliminaries}%}

First, we recall some definitions and relative properties of extriangulated category from \cite{NP}. Throughout this paper, we assume that $\mathcal{C}$ is an additive category.

\bg{Def}$\mathrm{\cite[Definition~2.1]{NP}}$\label{}
Assume that that C is equipped with an additive bifunctor $\mathbb{E}$: $C^{op}\times C\longrightarrow \mathrm{Ab}$, where $\mathrm{Ab}$ is the category of Abelian groups.
For any pair of objects $A,C\in\mathcal{C}$, an element $\delta \in \mathbb{E}(C,A)$ is
called an $\mathbb{E}$-extension. Thus formally, an $\mathbb{E}$-extension is a triplet $(A,\delta,C)$. If $\delta=0$, we called the split $\mathbb{E}$-extension.
\ed{Def}

For any morphism $a\in\mathcal{C}(A,A')$ and $c\in\mathcal{C}(C',C)$, we have $\mathbb{E}$-extensions
$\mathbb{E}(C,a)(\delta)\in\mathbb{E}(C,A')$ and $E(c,A)(\delta)\in\mathbb{E}(C',A)$.
We abbreviately denote them by $a_{\ast}\delta$ and $c^{\ast}\delta$ respectively.
A morphism $(a,c)$: $\delta\longrightarrow\delta'$ of $\mathbb{E}$-extensions is a pair of morphisms $a\in\mathcal{C}(A,A')$ and
$c\in\mathcal{C}(C,C')$ in $\mathcal{C}$ satisfying the equality $a_{\ast}\delta=c^{\ast}\delta'$.

\bg{Def}$\mathrm{\cite[Definition~2.7]{NP}}$\label{}
Let $A,C\in \mathcal{C}$ be any pair of objects. Two sequences of morphisms $\xymatrix{A\ar[r]^{x} &B\ar[r]^{y} &C}$ and $\xymatrix{A\ar[r]^{x'} &B'\ar[r]^{y'} &C}$ are said to be
equivalent if there is an isomorphism $b\in\mathcal{C}(B,B')$ such that $bx=x'$ and $y'b=y$. We denote the equivalence class of $\xymatrix{A\ar[r]^{x} &B\ar[r]^{y} &C}$
by $[\xymatrix{A\ar[r]^{x} &B\ar[r]^{y} &C}]$.
\ed{Def}

\bg{Def}$\mathrm{\cite[Definition~2.8]{NP}}$\label{}
$(1)$ For any $A,C\in \mathcal{C}$, we denote as $0=[\xymatrix@C=0.5cm{A\ar[r]^{(1\\0)} &A\bigoplus C\ar[r]^{(0,1)} &C}]$.

$(2)$ For any two classes $[\xymatrix{A\ar[r]^{x} &B\ar[r]^{y} &C}]$ and $[\xymatrix{A'\ar[r]^{x'} &B'\ar[r]^{y'} &C'}]$,
we denote as
$$\xymatrix{
[A\ar[r]^{x} &B\ar[r]^{y} &C] \bigoplus [A'\ar[r]^{~~~x'} &B'\ar[r]^{y'} &C']=
[A\bigoplus A'\ar[r]^{~~~~x\oplus x'} &B\bigoplus B'\ar[r]^{y\oplus y'} &C\bigoplus C']
}$$.
\ed{Def}

\bg{Def}$\mathrm{\cite[Definition~2.9]{NP}}$\label{}
Let $\mathfrak{s}$ be a correspondence which associates an equivalence
class $\mathfrak{s}(\delta)=[\xymatrix{A\ar[r]^{x} &B\ar[r]^{y} &C}]$ to any $\mathbb{E}$-extension $\delta\in\mathbb{E}(C,A)$.
This $\mathfrak{s}$ is called a realization of $\mathbb{E}$, if it satisfies the following condition $(\ast)$. In this
case, we say that sequence $[\xymatrix{A\ar[r]^{x} &B\ar[r]^{y} &C}]$ realizes $\delta$, whenever it satisfies
$\mathfrak{s}(\delta)=[\xymatrix{A\ar[r]^{x} &B\ar[r]^{y} &C}]$.
\ed{Def}

$(\ast)$ Let $\delta \in\mathbb{E}(C,A)$ and $\delta'\in\mathbb{E}(C',A')$ be any pair of $\mathbb{E}$-extensions, with
$\mathfrak{s}(\delta)=[\xymatrix@C=0.5cm{A\ar[r]^{x} &B\ar[r]^{y} &C}]$ and $\mathfrak{s}(\delta')=[\xymatrix{A'\ar[r]^{x'} &B'\ar[r]^{y'} &C'}]$.
Then, for any morphism $(a,c)$: $\delta\longrightarrow \delta'$, there exists $b\in\mathcal{C}(B,B')$ such that $bx=x'a$ and $y'b=cy$.

\bg{Def}$\mathrm{\cite[Definition~2.10]{NP}}$\label{}
A realization $\mathfrak{s}$ of $\mathbb{E}$ is said to be additive, if it satisfies the following conditions.

$(1)$ For any $A,C\in \mathcal{C}$, the split $\mathbb{E}$-extension $0\in\mathbb{E}(C,A)$ satisfies $\mathfrak{s}(0) = 0$.

$(2)$ For any pair of $\mathbb{E}$-extensions $\delta$ and $\delta'$,
$\mathfrak{s}(\delta\oplus\delta')=\mathfrak{s}(\delta)\oplus\mathfrak{s}(\delta')$ holds.
\ed{Def}

\bg{Def}$\mathrm{\cite[Definition~2.12]{NP}}$\label{etri}
We call the triplet $(\mathcal{C},\mathbb{E},\mathfrak{s})$ (simply, $\mathcal{ C}$) an extriangulated category if it satisfies the following conditions:

$\mathrm{(ET1)}$ $\mathbb{E}$: $C^{op}\times C\longrightarrow \mathrm{Ab}$ is an additive functor.

$\mathrm{(ET2)}$ $\mathfrak{s}$ is an additive realization of $\mathbb{E}$.

$\mathrm{(ET3)}$ Let $\delta\in\mathbb{E}(C,A)$ and $\delta'\in\mathbb{E}(C',A')$ be any pair of $\mathbb{E}$-extensions,
realized as $\mathfrak{s}(\delta)=[\xymatrix{A\ar[r]^{x} &B\ar[r]^{y} &C}]$ and $\mathfrak{s}(\delta')=[\xymatrix{A'\ar[r]^{x'} &B'\ar[r]^{y'} &C'}]$.
For any commutative square in $\mathcal{C}$
$$\xymatrix{
A\ar[r]^{x}\ar[d]^{a}&B\ar[r]^{y}\ar[d]^{b}&C\\
A'\ar[r]^{x'} &B'\ar[r]^{y'} &C'
}$$
there exists a morphism $(a,c)$: $\delta\longrightarrow\delta'$ which is realized by $(a,b,c)$.

$(\mathrm{ET3})$$^{\mathrm{op}}$  Dual of $\mathrm{(ET3)}$.

$(\mathrm{ET4})$ Let $(A,\delta,D)$ and $(B,\delta',F)$ be $\mathbb{E}$-extensions respectively realized by $\xymatrix{A\ar[r]^{f} &B\ar[r]^{f'} &D}$
and $\xymatrix{B\ar[r]^{g} &C\ar[r]^{g'} &F}$. Then there exist an object $E\in\mathcal{C}$, a commutative diagram
$$\xymatrix{
A\ar[r]^{f}\ar@{=}[d] &B\ar[r]^{f'}\ar[d]^{g}&D\ar[d]^{d}\\
A\ar[r]^{h} &C\ar[r]^{h'}\ar[d]^{g'} &E\ar[d]^{e}\\
&F\ar@{=}[r]&F
}$$
in $\mathcal{C}$, and an $\mathbb{E}$-extension $\delta''\in\mathbb{E}(E,A)$ realized by $\xymatrix{A\ar[r]^{h} &C\ar[r]^{h'} &E}$,
which satisfy the following compatibilities.

$(i)$ $\xymatrix{D\ar[r]^{d} &E\ar[r]^{e} &F}$ realizes $\mathbb{E}(F,f')(\delta')$,

$(ii)$ $\mathbb{E}(d,A)(\delta'') = \delta$,

$(iii)$ $\mathbb{E}(E,f)(\delta'') = \mathbb{E}(e,B)(\delta')$.

By (iii), $(f,e)$: $\delta''\longrightarrow\delta$ is a morphism of $\mathbb{E}$-extensions, realized by
$$(f,id_{C} ,e):~[\xymatrix{A\ar[r]^{h} &C\ar[r]^{h'} &E}]\longrightarrow [\xymatrix{B\ar[r]^{g} &C\ar[r]^{g'} &F}]$$.

$\mathrm{(ET4)}$$^{\mathrm{op}}$ Dual of $\mathrm{(ET4)}$.
\ed{Def}

In Definition \ref{etri}, if the sequence $\xymatrix{A\ar[r]^{x} &B\ar[r]^{y} &C}$ realizes $\delta\in\mathbb{E}(C,A)$,
we call the pair $(\xymatrix{A\ar[r]^{x} &B\ar[r]^{y} &C},\delta)$ an $\mathbb{E}$-triangle, and write it in the following way,
$$\xymatrix{A\ar[r]^{x} &B\ar[r]^{y} &C\ar@{.>}[r]^{\delta}&}.$$
And the sequence $\xymatrix{A\ar[r]^{x} &B\ar[r]^{y} &C}$ is called a conflation, $x$ is called an inflation and $y$ is called
a deflation.

\bg{Def}$\mathrm{\cite[Corollary~3.12]{NP}}$\label{ }
Let $\mathcal{C}$ be an extriangulated category. For any $\mathbb{E}$-triangle $\xymatrix{A\ar[r]  &B\ar[r] &C\ar@{.>}[r]^{\delta}&}$,
the following sequences of natural transformations are exact.
$$\xymatrix{\mathcal{C}(C,-)\ar[r] &\mathcal{C}(B,-)\ar[r] &\mathcal{C}(A,-)\ar[r]
&\mathbb{E}(C,-)\ar[r]&\mathbb{E}(B,-)\ar[r]&\mathbb{E}(A,-)},$$
$$\xymatrix{\mathcal{C}(-,A)\ar[r] &\mathcal{C}(-,B)\ar[r]&\mathcal{C}(-,C)\ar[r]
&\mathbb{E}(-,A)\ar[r]&\mathbb{E}(-,B)\ar[r]&\mathbb{E}(-,C)}.$$
\ed{Def}

In the rest of the section, we always assume that the additive category $\mathcal{ C}$ is an extriangulated category satisfying 
the following conditions, see \cite[Condition~5.8]{NP}.

\bg{Con}$\mathrm{(WIC)}$\label{ }
Let $\mathcal{C}$ be an extriangulated category, $f$ and $g$ be any composable pair of morphisms in $\mathcal{C}$.
Consider the following conditions.

$(1)$ If $gf$ is an inflation, then so is $f$.

$(2)$ If $gf$ is an deflation, then so is $g$.
\ed{Con}

\section{Gorenstein syzygy objects}
%---------------------------------------------------------------------------
\vskip 10pt
We first call some results in \cite{HZZ}. Here, we suppose that the proper classes of $\mathbb{E}$-triangles \cite{HZZ} is
all $\mathbb{E}$-triangles of $\mathcal{ C}$. We mainly give the definition of Gorenstein syzygy objects in extriangulated categories and
get a its characterization, see the theorem \ref{th1}.
And we can obtain a equivalent characterization about the Gorenstein projective dimension, see the theorem \ref{th2}.

\bg{Def}$\mathrm{\cite[Definition~4.1]{HZZ}}$\label{}
An object $P\in \mathcal{C}$ is called projective if for any  $\mathbb{E}$-triangle
$\xymatrix@C=0.5cm{A\ar[r] &B\ar[r] &C\ar@{.>}[r]^{\delta}&}$, the induced sequence of of abelian group
$0\longrightarrow\mathcal{C}(P,~A)\longrightarrow\mathcal{C}(P,~B)\longrightarrow\mathcal{C}(P,~C)\longrightarrow0$
is exact in $\mathrm{Ab}$. Dually, we can define injective objects.
\ed{Def}

We denote $\mathcal{P(C)}$ (resp. $\mathcal{I(C)}$) the class of projective (resp. injective) objects of $\mathcal{C}$.
An extrianglated category $(\mathcal{C},~\mathbb{E},~\mathfrak{s})$ is said to have enough projective (resp. enough injective) if for any object $A$
there is an $\mathbb{E}$-triangle $\xymatrix{B\ar[r] &P\ar[r] &A\ar@{.>}[r]&}$
(resp. $\xymatrix{A\ar[r] &I\ar[r] &B\ar@{.>}[r]&}$) with $P\in \mathcal{P(C)}$ (resp. $I\in \mathcal{I(C)}$).

\bg{Def}$\mathrm{\cite[Definition~4.4]{HZZ}}$\label{}
An unbounded complex $\mathbf{X}$ is called exact complex if $\mathbf{X}$ is a sequence
$$\xymatrix{
\cdots\ar[r]&X_{1}\ar[r]^{d_{1}} &X_{0}\ar[r]^{d_{0}} &X_{-1}\ar[r]^{d_{-1}} &\cdots\
}$$
in $\mathcal{C}$ such that there exists an $\mathbb{E}$-triangle
$\xymatrix{
K_{n+1}\ar[r]^{g_{n}} &X_{n}\ar[r]^{f_{n}} &K_{n}\ar@{.>}[r]^{\delta_{n}}&
}$
for any integer $n$ and $d_{n}=g_{n-1}f_{n}$.
\ed{Def}

\bg{Def}$\mathrm{\cite[Definition~4.5]{HZZ}}$\label{}
Let $\mathcal{W}$ be a class of objects in $\mathcal{C}$. An $\mathbb{E}$-triangle $\xymatrix@C=0.5cm{A\ar[r] &B\ar[r]&C\ar@{.>}[r]&}$ is
said to be $\mathcal{C}(-,~\mathcal{W})$-exact if for any $W\in \mathcal{W}$, the induced sequence of abelian group
$0\longrightarrow\mathcal{C}(C,~W)\longrightarrow\mathcal{C}(B,~W)\longrightarrow\mathcal{C}(A,~W)\longrightarrow0$ is exact in $\mathrm{Ab}$.
\ed{Def}

\bg{Def}$\mathrm{\cite[Definition~4.6]{HZZ}}$\label{}
Let $\mathcal{W}$ be a class of objects in $\mathcal{C}$. An complex $\mathbf{X}$ is said to be $\mathcal{C}(-,~\mathcal{W})$-exact
if it is an exact complex
$$\xymatrix{
\cdots\ar[r]&X_{1}\ar[r]^{d_{1}} &X_{0}\ar[r]^{d_{0}} &X_{-1}\ar[r]^{d_{-1}} &\cdots\
}$$
such that there exists an $\mathcal{C}(-,~\mathcal{W})$-exact $\mathbb{E}$-triangle
$\xymatrix{K_{n+1}\ar[r]^{g_{n}} &X_{n}\ar[r]^{f_{n}} &K_{n}\ar@{.>}[r]^{\delta_{n}}&}$
and $d_{n}=g_{n-1}f_{n}$ for any integer $n$.
\ed{Def}

Dually, we can define $\mathcal{C}(\mathcal{W},-)$-exact $\mathbb{E}$-triangles and $\mathcal{C}(\mathcal{W},-)$-exact complexes.
An exact complex $\mathbf{X}$ is called complete $\mathcal{W}$-exact if it is both $\mathcal{C}(-,~\mathcal{W})$-exact and $\mathcal{C}(\mathcal{W},-)$-exact.

\bg{Def}$\mathrm{\cite[Definition~4.7]{HZZ}}$\label{}
A complete projective resolution is a complete $\mathcal{P(C)}$-exact complex
$$\xymatrix{
\mathbf{P}:~\cdots\ar[r]&P_{1}\ar[r]^{d_{1}} &P_{0}\ar[r]^{d_{0}} &P_{-1}\ar[r]^{d_{-1}} &\cdots\
}$$
in $\mathcal{C}$ such that $P_{i}$ is projective for any integer $i$. Dually, we can define complete injective resolution.
\ed{Def}

\bg{Def}$\mathrm{\cite[Definition~4.8]{HZZ}}$\label{}
Let $\mathbf{P}$ be a complete projective resolution in $\mathcal{C}$. So there exists a $\mathcal{C}(-,\mathcal{P(C)})$-exact
$\mathbb{E}$-triangle $\xymatrix{K_{n+1}\ar[r]^{g_{n}} &P_{n}\ar[r]^{f_{n}} &K_{n}\ar@{.>}[r]^{\delta_{n}}&}$ for any integer $n$.
The objects $K_{n}$ are called Gorenstein projective for any $n$. Dually, we can define complete injective resolution and Gorenstein
injective objects. We denote $\mathcal{GP(C)}$ (resp. $\mathcal{GI(C)}$) the class of Gorenstein projective (resp. Gorenstein injective)
objects in $\mathcal{C}$.
\ed{Def}

\bg{Def}\label{GSYZ}
Let $n$ be a positive integer. If there is an exact complex
$$\xymatrix{
0\ar[r]&N\ar[r] &G_{n-1}\ar[r] &\cdots \ar[r]&G_{1}\ar[r]&G_{0}\ar[r]&M\ar[r]&0
}$$
with $G_{i}\in \mathcal{GP(C)}$ for any $i$, then $N$ is called a Gorenstein $n$-syzygy object (of $M$).
\ed{Def}

In the definition above, if $G_{i}\in \mathcal{P(C)}$ for any $i$, then the object $N$ is called $n$-syzygy (of $M$).

\bg{Lem}\label{gp}
Let $N$ be a Gorenstein $2$-syzygy object of $M$ in $\mathcal{C}$. Then we can obtain the following two exact complexes
$$\xymatrix{
0\ar[r]&N\ar[r] &H_{1}\ar[r]&P_{1}\ar[r]&M\ar[r]&0
}$$
and
$$\xymatrix{
0\ar[r]&N\ar[r] &P_{2}\ar[r]&H_{2}\ar[r]&M\ar[r]&0
}$$
with $P_{i}\in \mathcal{P(C)}$ and $H_{i}\in \mathcal{GP(C)}$ for $i=1,2$.
\ed{Lem}

\Pf. Since $N$ is Gorenstein $2$-syzygy, there is an exact complex
$$\xymatrix{0\ar[r]&N\ar[r] &G_{1}\ar[r]&G_{0}\ar[r]&M\ar[r]&0}$$
with $G_{i}\in \mathcal{GP(C)}$ for $i=1,2$.
And then there exist two $\mathbb{E}$-triangles $\xymatrix{N\ar[r] &G_{1}\ar[r] &K\ar@{.>}[r] &}$ and
$\xymatrix{K\ar[r] &G_{0}\ar[r] &M\ar@{.>}[r] &}$. Because $G_{0}\in \mathcal{GP(C)}$, we have an $\mathbb{E}$-triangles
$\xymatrix@C=0.5cm{G'_{0}\ar[r] &P_{1}\ar[r] &G_{0}\ar@{.>}[r] &}$. By (ET4)$^{\mathrm{op}}$, we have the following commutative diagrams.
$$\xymatrix{
G'_{0}\ar@{=}[r]\ar[d]&G'_{0}\ar[d]&\\
L\ar[r]^{y}\ar[d]&P_{1}\ar[r]\ar[d]&M\ar@{.>}[r]\ar@{=}[d]&\\
K\ar[r]\ar@{.>}[d]&G_{0}\ar[r]\ar@{.>}[d]&M\ar@{.>}[r]&\\
&&&
}$$
and
$$\xymatrix{
&G'_{0}\ar@{=}[r]\ar[d]&G'_{0}\ar[d]\\
N\ar[r]\ar@{=}[d]&H_{1}\ar[r]^{x}\ar[d]&L\ar@{.>}[r]\ar[d]&\\
N\ar[r]&G_{1}\ar[r]\ar@{.>}[d]&K\ar@{.>}[r]\ar@{.>}[d]&\\
&&
}$$
From the above diagrams, we can obtain the following exact complex.
$$\xymatrix{
0\ar[r]&N\ar[r] &H_{1}\ar[r]^{yx}&P_{1}\ar[r]&M\ar[r]&0
}$$
In the second diagram, we have that $H_{1}\in \mathcal{GP(C)}$ since $\mathcal{GP(C)}$ is closed under extensions \cite[Theorem~4.16]{HZZ}.
Dually, we can get the other exact complex.
\ \hfill $\Box$

\mskip\

Next, we give the first result in this section.

\bg{Th}\label{th1}
Let $n$ be a positive integer and $N\in\mathcal{ C}$. Then the following are equivalent.

$(1)$ There is an object $M\in C$ such that N is a Gorenstein $n$-syzygy object of $M$. i.e., there is an exact complex
$$\xymatrix{
0\ar[r]&N\ar[r] &G_{n-1}\ar[r]^{d_{n-1}} &\cdots \ar[r]&G_{1}\ar[r]^{d_{1}}&G_{0}\ar[r]^{d_{0}}&M\ar[r]&0
}$$
with $G_{i}\in \mathcal{GP(C)}$ for all $i$.

$(2)$ There is a deflation $L\longrightarrow G$ in $\mathcal{C}$ such that
$N$ is an $n$-syzygy object of $L$ and $G\in\mathcal{GP(C)}$.
\ed{Th}

\Pf. $(2)\Rightarrow (1)$ Note that $\mathcal{P(C)}\subseteq\mathcal{GP(C)}$. Taking $M=L$, (1) holds.

$(1)\Rightarrow (2)$ We use induction on $n$. When $n=1$, there is an E-triangle $\xymatrix{G_{0}\ar[r] &P_{0}\ar[r] &G\ar@{.>}[r] &}$ since $G_{0}\in\mathcal{GP(C)}$.
We consider the following $\mathbb{E}$-triangles commutative diagram.
$$\xymatrix{
N\ar[r]\ar@{=}[d]&G_{0}\ar[r]\ar[d]&M\ar[d]\ar@{.>}[r]&\\
N\ar[r]&P_{0}\ar[r]\ar[d]&L\ar@{.>}[r]\ar[d]&\\
&G\ar@{=}[r]\ar@{.>}[d]&G\ar@{.>}[d]&\\
&&&
}$$
From the second row and third column in diagram above, we can know that (1) holds.

Assume that the result holds for the case $n-1$, where $n\geq2$. Since $N$ is Gorenstein $n$-syzygy, there is an exact complex
$$(\ast1):\xymatrix{
0\ar[r]&K_{n-2}\ar[r]^{g_{n-3}} &G_{n-3}\ar[r]^{d_{n-3}} &\cdots \ar[r]&G_{1}\ar[r]^{d_{1}}&G_{0}\ar[r]^{d_{0}}&M\ar[r]&0
}$$
and an $\mathbb{E}$-triangle
$$\xymatrix{
0\ar[r]&N\ar[r] &G_{n-1}\ar[r] &G_{n-2}\ar[r]^{f_{n-2}} &K_{n-2}\ar[r]&0
}$$
where $d_{n-2}=g_{n-3}f_{n-2}$. By the lemma \ref{gp}, we can obtain the following exact complex
$$(\ast2):\xymatrix{
0\ar[r]&N\ar[r] &P_{n-1}\ar[r]^{d'_{n-1}} &G'_{n-2}\ar[r]^{f'_{n-2}} &K_{n-2}\ar[r]&0
}$$
with $P_{n-1}\in \mathcal{P(C)}$ and $G'_{n-2}\in \mathcal{GP(C)}$. From the two exact complexes $(\ast1)$ and $(\ast2)$, we obtain
the following exact complex.
$$(\ast3):\xymatrix@C=0.5cm{
0\ar[r]&N\ar[r] &P_{n-1}\ar[r]^{d'_{n-1}} \ar[r]&G'_{n-2}\ar[r]^{d'_{n-2}}&G_{n-3}\ar[r]^{d_{n-3}}&\cdots \ar[r]&G_{1}\ar[r]^{d_{1}}&G_{0}\ar[r]^{d_{0}}&M\ar[r]&0
}$$
with $d'_{n-2}=g_{n-3}f'_{n-2}$. From $(\ast3)$, we can obtain the following exact complex
$$\xymatrix{
0\ar[r]&K'_{n-1}\ar[r]^{g'_{n-1}} &G'_{n-2}\ar[r]^{d'_{n-2}} &G_{n-3}\ar[r]^{d_{n-3}} &\cdots \ar[r]&G_{0}\ar[r]^{d_{0}}&M\ar[r]&0
}$$
and an $\mathbb{E}$-triangle
$$(\ast4):\xymatrix{
N\ar[r] &P_{n-1}\ar[r]^{f'_{n-1}} &K'_{n-1}\ar@{.>}[r]&
}$$
where $d'_{n-1}=g'_{n-1}f'_{n-1}$. By the inductive hypothesis, we can know that $K'_{n-1}$ is an $(n-1)$-syzygy object of $L$. i.e., we have the following exact complex.
$$(\ast5):\xymatrix{
0\ar[r]&K'_{n-1}\ar[r] &P_{n-2}\ar[r]&\cdots \ar[r]&P_{1}\ar[r]&P_{0}\ar[r]&L\ar[r]&0
}$$
From $(\ast4)$ and $(\ast5)$, we know that $N$ is an $n$-syzygy object of $L$.
Existence of $G\in \mathcal{GP(C)}$, similar to the case $n=1$.
\ \hfill $\Box$

\mskip\

The Gorenstein projective dimension $\mathrm{Gpd}M$ of an object $M \in \mathcal{C}$ is defined inductively.
If $M\in \mathcal{GP(C)}$ then define $\mathrm{Gpd}M=0$. Next by induction, for an integer $n>0$, put $\mathrm{Gpd}M\leq n$ if there is an
$\mathbb{E}$-triangle $\xymatrix{N\ar[r] &G\ar[r]&M\ar@{.>}[r]&}$ with $G\in \mathcal{GP(C)}$ and $\mathrm{Gpd}N\leq n-1$.
Next, we will give a characterization of the Gorenstein projective dimension.

\bg{Th}\label{th2}
Let $M$ be an object in $\mathcal{C}$ and $n$ be a positive integer. Then $\mathrm{Gpd}M\leq n$ if and only if
for any integer $s$ satisfying $0\leq s\leq n$, there is an exact complex
$$\xymatrix{
0\ar[r]&X_{n}\ar[r] &X_{n-1}\ar[r] &\cdots \ar[r]&X_{1}\ar[r]&X_{0}\ar[r]&M\ar[r]&0
}$$
with $X_{s} \in \mathcal{GP(C)}$ and $X_{i}\in \mathcal{P(C)}$ for $i\neq s$.
\ed{Th}

\Pf. ($\Leftarrow$) Obviously.

($\Rightarrow$) We use induction on $n$. When $n=1$, taking $N=0$ in the lemma \ref{gp}, it is easy to see
the result holds.

Suppose that the result holds for the case $n-1$, where $n\geq2$. Since $\mathrm{Gpd}M\leq n$, then there is an exact complex
$$\xymatrix{
0\ar[r]&G_{n}\ar[r] &G_{n-1}\ar[r]^{d_{n-1}} &\cdots \ar[r]&G_{1}\ar[r]^{d_{1}}&G_{0}\ar[r]^{d_{0}}&M\ar[r]&0
}$$
with $G_{i} \in \mathcal{GP(C)}$ for any $i$. And then we have the following two exact complexes
$$(\sharp1):\xymatrix{
0\ar[r]&G_{n}\ar[r]&G_{n-1}\ar[r]^{d_{n-1}} &\cdots \ar[r]&G_{2}\ar[r]^{f_{2}}&K_{2}\ar[r]&0
}$$
and
$$(\sharp2):\xymatrix{
0\ar[r]&K_{2}\ar[r]^{g_{1}} &G_{1}\ar[r] &G_{0}\ar[r]&M\ar[r]&0
}$$
where $d_{2}=g_{1}f_{2}$. Since $K_{2}$ is Gorenstein $2$-syzygy from $(\sharp2)$, by the lemma \ref{gp}, we can obtain the following exact complex
$$(\sharp3):\xymatrix{
0\ar[r]&K_{2}\ar[r]^{g'_{1}} &G\ar[r] &P\ar[r]&M\ar[r]&0
}$$
with $G\in\mathcal{GP(C)}$ and $P\in\mathcal{P(C)}$. Thus, from $(\sharp1)$ and $(\sharp3)$, we have the following exact complexes
$$(\sharp4):\xymatrix{
0\ar[r]&G_{n}\ar[r]&G_{n-1}\ar[r]^{d_{n-1}} &\cdots \ar[r]&G_{2}\ar[r]^{d'_{2}}&G\ar[r]^{d'_{1}}&P\ar[r]&M\ar[r]&0
}$$
where $d'_{2}=g'_{1}f_{2}$. From $(\sharp4)$, we have the following exact complexes
$$\xymatrix{
0\ar[r]&G_{n}\ar[r]&G_{n-1}\ar[r]^{d_{n-1}} &\cdots \ar[r]&G_{2}\ar[r] &G\ar[r]^{f'_{1}}&K'_{1}\ar[r]&0
}$$
and an $\mathbb{E}$-triangle
$$(\sharp5):\xymatrix{
K'_{1}\ar[r]^{g'_{0}} &P\ar[r] &M\ar@{.>}[r]&
}$$
where $d'_{1}=g'_{0}f'_{1}$. Note that $\mathrm{Gpd}K'_{1}\leq n-1$. By the inductive hypothesis and $(\sharp5)$, the result holds for $s\neq0$.

Next, we only need to prove the result for $s=0$.
We have the following an exact complexes
$$(\natural1):\xymatrix{
0\ar[r]&G_{n}\ar[r]&G_{n-1}\ar[r]^{d_{n-1}} &\cdots \ar[r]&G_{1}\ar[r]^{f_{1}}&K_{1}\ar[r]&0
}$$
and an $\mathbb{E}$-triangle
$$(\natural2):\xymatrix{
0\ar[r]&K_{1}\ar[r]^{g_{0}} &G_{0}\ar[r] &M\ar[r]&0
}$$
where $d_{1}=g_{0}f_{1}$. By the inductive hypothesis, we have the following an exact complexes
$$(\natural3):\xymatrix{
0\ar[r]&P_{n}\ar[r]&P_{n-1}\ar[r] &\cdots \ar[r]&P_{2}\ar[r]^{\sigma_{2}}&G'_{1}\ar[r]^{f''_{1}} &K_{1}\ar[r]&0
}$$
with $P_{i}\in \mathcal{P(C)}$ and $G'_{1}\in \mathcal{GP(C)}$. From $(\natural2)$ and $(\natural3)$, we have the following an exact complexes
$$(\natural4):\xymatrix{
0\ar[r]&P_{n}\ar[r]&P_{n-1}\ar[r] &\cdots \ar[r]&P_{2}\ar[r]^{\sigma_{2}}&G'_{1}\ar[r]^{g_{0}f''_{1}} &G_{0}\ar[r] &M\ar[r]&0
}$$
Consider the following two exact complexes
$$(\natural5):\xymatrix{
0\ar[r]&P_{n}\ar[r]&P_{n-1}\ar[r] &\cdots \ar[r]&P_{2}\ar[r]^{f''_{2}}&K''_{2}\ar[r]&0
}$$
and
$$(\natural6):\xymatrix{
0\ar[r]&K''_{2}\ar[r]^{g''_{1}}&G'_{1}\ar[r]^{g_{0}f''_{1}} &G_{0}\ar[r] &M\ar[r]&0
}$$
where $\sigma_{2}=g''_{1}f''_{2}$. Since $K''_{2}$ is Gorenstein $2$-syzygy from $(\natural6)$, by the lemma \ref{gp}, we can obtain the following exact complex
$$(\natural7):\xymatrix{
0\ar[r]&K''_{2}\ar[r] &P''\ar[r] &G''\ar[r]&M\ar[r]&0
}$$
with $G''\in\mathcal{GP(C)}$ and $P''\in\mathcal{P(C)}$. From $(\natural5)$ and $(\natural7)$, the result holds.
\ \hfill $\Box$

\section{Recollement}

Let the functor $F$: ${ \mathcal{A}}\longrightarrow {\mathcal{B}}$ and the functor $G$: $\mathcal{B}\longrightarrow \mathcal{A}$ are additive,
where both $\mathcal{A}$ and $\mathcal{B}$ are extriangulated categories. we said $(F$, $G)$ to be an adjoint pair,
if there is a isomorphism $\sigma_{X,Y}$: ${\mathcal{B}}(FX, Y)\cong{\mathcal{A}}(X, GY)$ for any $X\in \mathcal{A}$ and $Y\in \mathcal{B}$.
In the case, if there is a functor $H$: $\mathcal{A}\longrightarrow \mathcal{B}$ such that the pair $(G,~H)$ is also an adjoint pair,
then we call $(F,~G,~H)$ being adjoint triple.

Let (${\mathcal{A}}$, $\mathbb{E}_{\mathcal{A}}$, $\mathfrak{s}_{\mathcal{A}}$) and ($\mathcal{B}$, $\mathbb{E}_{\mathcal{B}}$, $\mathfrak{s}_{\mathcal{B}}$)
be extriangulated categories. We say an additive covariant functor $F$: $\mathcal{A}\longrightarrow \mathcal{B}$ is an exact functor
if it preserves $\mathbb{E}$-triangles. First, we give the definition of recollement of extriangulated categories.

\bg{Def}$\mathrm{\cite[Definition~3.1]{WWZ}}$\label{recollo}%
Let $\mathcal{A}$, $\mathcal{A'}$, $\mathcal{A''}$ be three extriangulated categories. A recollement of $\mathcal{A}$ relative $\mathcal{A'}$ and
$\mathcal{A''}$, denoted by $R(\mathcal{A'}, ~\mathcal{A}, ~\mathcal{A''})$, is a diagram
$$\xymatrix{
{\mathcal{A}^\prime} \ar[rrr]|{\ i_{\ast}}&&& {\mathcal{A}}
\ar@/^/@<2ex>[lll]|{i^{!}}\ar[rrr]|{\
j^{\ast}}\ar@/_/@<-2ex>[lll]|{i^{\ast}} &&& {\mathcal{A''}}
\ar@/_/@<-2ex>[lll]|{j_{!}}\ar@/^/@<2ex>[lll]|{j_{\ast}} }$$
given by two exact functors $i_{\ast}$, $j^{\ast}$, two right exact functors $i^{\ast}$, $j_{!}$ and two left exact functors $i^{!}$, $j_{\ast}$,
which satisfies the conditions below:

$(R1)$ ($i^{\ast}$, $i_{\ast}$, $i^{!}$) and ($j_{!}$, $j^{\ast}$, $j_{\ast}$) are adjoint triples;

$(R2)$ The functors $i_{\ast}$, $j_{!}$, and $j_{\ast}$ are fully faithful;

$(R3)$ $\mathrm{Im}i_{\ast}$ = $\mathrm{Ker}j^{\ast}$;

$(R4)$ For any $X\in\mathcal{A}$, there is a left exact $\mathbb{E}_{\mathcal{A}}$-triangle sequence
$$\xymatrix{
i_{\ast}i^{!}X\ar[r]^{\theta_{X}}&X\ar[r]^{\vartheta_{X}} &j_{\ast}j^{\ast}X\ar[r] &i_{\ast}A
}$$
with $A\in\mathcal{A'}$, where $\theta_{X}$ and $\vartheta_{X}$ are given by the adjunction morphisms.

$(R5)$ For any $X\in\mathcal{A}$, there is a right exact $\mathbb{E}_{\mathcal{A}}$-triangle sequence
$$\xymatrix{
i_{\ast}A'\ar[r]&j_{!}j^{\ast}X\ar[r]^{\alpha_{X}} &X\ar[r]^{\beta X} &i_{\ast}i^{\ast}X
}$$
with $A'\in\mathcal{A'}$, where $\alpha_{X}$ and $\beta_{X}$ are given by the adjunction morphisms.
\ed{Def}

Here, for the definition of the left (right) exact sequence (resp., functor), please refer to \cite{WWZ} for details.
Since these definitions are not used in this article, their definitions are not given.

\bg{Rem}\label{}%
$(1)$ If the categories $\mathcal{A}$, $\mathcal{A'}$ and $\mathcal{A''}$ are Abelian, then Definition \ref{recollo} coincides with
the definition of recollement of Abelian categories (cf. $\mathrm{\cite{MH,PC,FP}}$)

$(2)$ If the categories $\mathcal{A}$, $\mathcal{A'}$ and $\mathcal{A''}$ are triangulated categories, then Definition \ref{recollo} coincides with
the definition of recollement of triangulated categories (cf. $\mathrm{\cite{BBD}}$)
\ed{Rem}

Next, we collect some properties of recollements, which is very useful, see \cite{WWZ}.

\bg{Pro}\label{property}%

Let $R(\mathcal{A'}, ~\mathcal{A}, ~\mathcal{A''})$ be a recollement of extriangulated categories. Then we have the following properties.

$(1)$ $i^{\ast}j_{!}=0$ and $i^{!}j_{\ast}=0$;

$(2)$  $i^{\ast}$ preserves projective objects and $i^{!}$  preserves injective objects;

$(2')$  $j_{!}$ preserves projective objects and $j_{\ast}$  preserves injective objects;

$(3)$ These natural transformations $i^{\ast}i_{\ast}\longrightarrow \mathrm{Id}_{\mathcal{A'}}$, $\mathrm{Id}_{\mathcal{A'}}\longrightarrow i^{!}i_{\ast}$,
$j^{\ast}j_{\ast}\longrightarrow \mathrm{Id}_{\mathcal{A''}}$, $\mathrm{Id}_{\mathcal{A''}}\longrightarrow j^{\ast}j_{!}$ are natural isomorphisms.

$(4)$ If $i^{\ast}$ is exact, then $j_{!}$ is exact;

$(4')$ If $i^{!}$ is exact, then $j_{\ast}$ is exact.

$(5)$ If $\mathcal{A}$ has enough projective objects, then $\mathcal{A'}$ has enough projective objects and $\mathcal{P}(\mathcal{A'})=add(i^{\ast}(\mathcal{P}(\mathcal{A})))$;
If $\mathcal{A}$ has enough injective objects, then $\mathcal{A'}$ has enough injective objects and $\mathcal{I}(\mathcal{A'})=add(i^{!}(\mathcal{I}(\mathcal{A})))$.

$(6)$ If $\mathcal{A}$ has enough projective objects and $j_{\ast}$ is exact, then $\mathcal{A}''$ has enough projective objects and
$\mathcal{P}(\mathcal{A}'')=add(j^{\ast}(\mathcal{P}(\mathcal{A})))$;
If $\mathcal{A}$ has enough injective objects and $j_{!}$ is exact, then $\mathcal{A}''$ has enough injective objects and
$\mathcal{I}(\mathcal{A}'')=add(j^{\ast}(\mathcal{I}(\mathcal{A})))$.

\ed{Pro}

\mskip\

In order to prove the main result of this section, we need the following two propositions.

\bg{Pro}\label{Pro1}
In a recollement $R(\mathcal{A}', ~\mathcal{A}, ~\mathcal{A}'')$, let $X$ be a Gorenstein projective
object in $\mathcal{A}$.

$(1)$ If $j_{\ast}$ is an exact functor, then $j^{\ast}X\in \mathcal{GP(A'')}$.

$(2)$ If $i^{\ast}$ is an exact functor, then $i^{\ast}X\in \mathcal{GP(A')}$.
\ed{Pro}

\Pf. (1) Since $X$ is a Gorenstein projective object in $\mathcal{A}$, there exists a $\mathcal{A}(-,\mathcal{P(A)})$-exact
$\mathbb{E}$-triangle $\xymatrix{X_{n+1}\ar[r] &P_{n}\ar[r] &X_{n}\ar@{.>}[r]&}$ for any integer $n$,
where each $P_{n}\in\mathcal{P(A)}$ and $X_{0}=X$. Note that $j^{\ast}j_{!}\longrightarrow \mathrm{Id}_{\mathcal{A''}}$ is natural isomorphism. We have the following exact commutative diagram.
$$\xymatrix{
\mathcal{A''}(j^{\ast}X_{n},~P'')\ar[r]\ar^{\cong}[d]&\mathcal{A''}(j^{\ast}P_{n},~P'')\ar[r]\ar^{\cong}[d]&\mathcal{A''}(j^{\ast}X_{n+1},~P'')\ar^{\cong}[d]\\
\mathcal{A''}(j^{\ast}X_{n},~j^{\ast}j_{!}P'')\ar[r]&\mathcal{A''}(j^{\ast}P_{n},~j^{\ast}j_{!}P'')\ar[r]&\mathcal{A''}(j^{\ast}X_{n+1},~j^{\ast}j_{!}P'')
}$$
where $P''\in\mathcal{GP(A'')}$. Note that $j_{!}P''\in\mathcal{GP(A)}$ by the proposition \ref{property}. So the sequence
$$\xymatrix{
0\ar[r]r&\mathcal{A''}(X_{n},~j_{!}P'')\ar[r]&\mathcal{A''}(P_{n},~j_{!}P'')\ar[r]&\mathcal{A''}(X_{n+1},~j_{!}P'')\ar[r]&0
}$$
is a short exact sequence. Thus the second row is a short exact sequence in above commutative diagram.
And then these $\mathbb{E}_{\mathcal{A''}}$-triangle $\xymatrix{j^{\ast}X_{n+1}\ar[r] &j^{\ast}P_{n}\ar[r] &j^{\ast}X_{n}\ar@{.>}[r]&}$ is $\mathcal{A''}(-,\mathcal{GP(A'')})$-exact,
where $j_{\ast}$ is an exact functor. By the proposition \ref{property}, $j^{\ast}P_{n}\in\mathcal{GP(A'')}$. Consequently, $j^{\ast}X\in \mathcal{GP(A'')}$.

(2) Similar to the discuss above, we have the following exact commutative diagram.
$$\xymatrix{
\mathcal{A'}(i^{\ast}X_{n},~P')\ar[r]\ar^{\cong}[d]&\mathcal{A'}(i^{\ast}P_{n},~P')\ar[r]\ar^{\cong}[d]&\mathcal{A'}(i^{\ast}X_{n+1},~P')\ar^{\cong}[d]\\
\mathcal{A}(X_{n},~i_{\ast}P')\ar[r]&\mathcal{A}(P_{n},~i_{\ast}P')\ar[r]&\mathcal{A}(X_{n+1},~i_{\ast}P')
}$$
where $P'\in\mathcal{GP(A')}$. It is easy to see that  the second row is a short exact sequence in above commutative diagram.
And then these $\mathbb{E}_{\mathcal{A'}}$-triangle $\xymatrix{i^{\ast}X_{n+1}\ar[r] &i^{\ast}P_{n}\ar[r] &i^{\ast}X_{n}\ar@{.>}[r]&}$ is $\mathcal{A'}(-,\mathcal{GP(A')})$-exact,
Note that $i^{\ast}P_{n}\in\mathcal{GP(A')}$ by the proposition \ref{property}. Consequently, $i^{\ast}X\in \mathcal{GP(A')}$.
\ \hfill $\Box$

\mskip\

We need the following lemma for the sake of the second proposition.

\bg{Lem}\label{lem0}
Let $\mathscr A$ and $\mathscr B$ be abelian categories, and let
$$L^{\bullet}=:\xymatrix{\cdots\ar[r]&L_{-1}\ar[r]^{d_{-1}}&L_{0}\ar[r]^{d_{0}}&L_{1}\ar[r]&\cdots}$$
be in ${\mathscr A}$. If the functor $F:$ $\mathscr A \longrightarrow \mathscr B$ is exact and faithful. Then $L^{\bullet}$ is exact if and only if $F(L^{\bullet})$ is exact.
\ed{Lem}

\Pf. The necessity is obvious since $F$ is exact.

Since $F(L^{\bullet})=:$ $\xymatrix{\cdots\ar[r]&FL_{-1}\ar[r]^{Fd_{-1}}&FL_{0}\ar[r]^{Fd_{0}}&FL_{1}\ar[r]&\cdots}$ is exact and $F$ is exact,
then the $i$-th homology  0=$\mathrm{H}_{i}(F(L^{\bullet}))=\mathrm{Ker} Fd_{i}/\mathrm{Im} Fd_{i-1}\cong F(\mathrm{Ker} d_{i}/\mathrm{Im} d_{i-1})$.
Note that an exact functor $F$ is faithful if and only if  it does not take any non-zero object to zero. So $\mathrm{Ker} d_{i}/\mathrm{Im} d_{i-1}=0$.
i.e., $L^{\bullet}$ is exact.
\ \hfill $\Box$

\bg{Pro}\label{Pro2}
In a recollement $R(\mathcal{A}', ~\mathcal{A}, ~\mathcal{A}'')$, let $X'$, $X''$ be a Gorenstein projective
object in $\mathcal{A}'$, $\mathcal{A}''$ respectively.

$(1)$ If $i^{!}$, $i^{\ast}$ are exact functors and $i^{\ast}$ is faithful, then $i_{\ast}X'\in \mathcal{GP(A)}$.

$(2)$ If $j_{!}$ and $j_{\ast}$ are exact functors, then $j_{!}X''\in \mathcal{GP(A)}$.
\ed{Pro}

\Pf. (1) Since $X'$ is a Gorenstein projective object in $\mathcal{A'}$, there exists a $\mathcal{A'}(-,\mathcal{P(A')})$-exact
$\mathbb{E}_{\mathcal{A'}}$-triangle $\xymatrix{X'_{n+1}\ar[r] &P'_{n}\ar[r] &X'_{n}\ar@{.>}[r]&}$ for any integer $n$,
where each $P'_{n}\in\mathcal{P(A')}$ and $X'_{0}=X'$. For any $P\in \mathcal{P(A)}$, applying the functor $i^{\ast}$ to the complex
$$(\ddag):\xymatrix{
0\ar[r]&\mathcal{A}(i_{\ast}X'_{n},P)\ar[r]&\mathcal{A}(i_{\ast}P'_{n},~P)\ar[r]&\mathcal{A}(i_{\ast}X'_{n+1},~P)\ar[r]&0
}$$
We have that the following complex is exact since $i^{\ast}P\in\mathcal{GP(A')}$.
$$\xymatrix{
0\ar[r]&\mathcal{A}(X'_{n},i^{\ast}P)\ar[r]&\mathcal{A}(P'_{n},~i^{\ast}P)\ar[r]&\mathcal{A}(X'_{n+1},~i^{\ast}P)\ar[r]&0
}$$
By the lemma \ref{lem0}, we have the sequence $(\ddag)$ is exact. By the proposition \ref{property}, $i_{\ast}P_{n}\in\mathcal{GP(A)}$.
Consequently, $i_{\ast}X'\in \mathcal{GP(A)}$.

(2) Since $X''$ is a Gorenstein projective object in $\mathcal{A'}$, there exists a $\mathcal{A''}(-,\mathcal{P(A'')})$-exact
$\mathbb{E}_{\mathcal{A''}}$-triangle $\xymatrix{X''_{n+1}\ar[r] &P''_{n}\ar[r] &X''_{n}\ar@{.>}[r]&}$ for any integer $n$,
where each $P''_{n}\in\mathcal{P(A'')}$ and $X''_{0}=X''$. We have the following exact commutative diagram.
$$\xymatrix{
\mathcal{A}(j_{!}X''_{n},~P)\ar[r]\ar^{\cong}[d]&\mathcal{A}(j_{!}P''_{n},~P)\ar[r]\ar^{\cong}[d]&\mathcal{A}(j_{!}X''_{n+1},~P)\ar^{\cong}[d]\\
\mathcal{A''}(X''_{n},~j^{\ast}P)\ar[r]&\mathcal{A''}(P''_{n},~j^{\ast}P)\ar[r]&\mathcal{A''}(X''_{n+1},~j^{\ast}P)
}$$
where $P\in\mathcal{GP(A)}$. Note that $j^{{\ast}}P\in\mathcal{GP(A'')}$.
So the second row is a short exact sequence in above commutative diagram.
And then these $\mathbb{E}_{\mathcal{A''}}$-triangle $\xymatrix{j^{\ast}X_{n+1}\ar[r] &j^{\ast}P_{n}\ar[r] &j^{\ast}X_{n}\ar@{.>}[r]&}$ is $\mathcal{A''}(-,\mathcal{GP(A'')})$-exact,
where $j_{\ast}$ is an exact functor. By the proposition \ref{property}, $j_{!}P''_{n}\in\mathcal{GP(A)}$. Consequently, $j_{!}X''\in \mathcal{GP(A)}$.
\ \hfill $\Box$

\mskip\

By the proposition \ref{Pro1} and \ref{Pro2}, we can get the main result of this section.

\bg{Th}\label{Th1}
In a recollement $R(\mathcal{A}', ~\mathcal{A}, ~\mathcal{A}'')$, let $X,~X',~X''$ be a Gorenstein n-syzygy
object in $\mathcal{A},~\mathcal{A'},~\mathcal{A''}$ respectively.

$(1)$ If $j_{\ast}$ is an exact functor, then $j^{\ast}X$ is Gorenstein n-syzygy in $\mathcal{A''}$.

$(2)$ If $i^{\ast}$ is an exact functor, then $i^{\ast}X$ is Gorenstein n-syzygy $\mathcal{A'}$.

$(3)$ If $i^{!}$, $i^{\ast}$ are exact functors and $i^{\ast}$ is faithful, then $i_{\ast}X'$ is Gorenstein n-syzygy $\mathcal{A}$.

$(4)$ If $j_{!}$ and $j_{\ast}$ are exact functors, then $j_{!}X''$ is Gorenstein n-syzygy $\mathcal{A}$.
\ed{Th}

{\small

}

\end{document}